# On Equivalence between Optimality Criteria and Projected Gradient Methods with Application to Topology Optimization Problem

SERGEY ANANIEV

Institute of Lightweight Structures and Conceptual Design
University of Stuttgart
Pfaffenwaldring 7
70569 Germany
Tel.: 49-711-685-6255
e-mail: sergey.ananiev@po.uni-stuttgart.de

**Abstract.** The paper demonstrates the equivalence between the optimality criteria (OC) method, initially proposed by Bendsøe & Kikuchi for topology optimization problem, and the projected gradient method. The equivalence is shown using Hestenes definition of Lagrange multipliers. Based on this development, an alternative formulation of the Karush-Kuhn-Tucker (KKT) condition is suggested. Such reformulation has some advantages, which will be also discussed in the paper. For verification purposes the modified algorithm is applied to the design of tension-only structures.

## 1. OC method for minimization of strain energy

Topology optimization, beginning with the pioneering work of Bendsøe & Kikuchi [1], where it was formulated as an optimal material distribution problem, has won a broad acceptance in industry and scientific world. For a recent overview of the successful applications of this technology the paper of Soto [2] can be consulted. Despite these dramatic developments, the original optimality criteria algorithm, which was used for numerical solution, remains without major changes. Even in the recent works in this field (Bendsøe & Sigmund [3], Maute [4]) it is still used as a "working horse" because of its simplicity and stability.

The interesting thing in this OC algorithm that in spite of the fact that it is derived from the Karush-Kuhn-Tucker optimality condition it uses the heuristic update rule, without any mathematical justification. This applies not only to the particular algorithm developed by Bendsøe & Kikuchi, but also to the generalized OC algorithm as presented by Venkayya [5] (see also Haftka [6] for overview).

To the author's knowledge, until now there was no critical analysis of the question *why* the OC methods work. This paper tries to find an answer to this question.

At first, the optimality criteria for minimum compliance design will derived. For simplicity reasons the isotropic material model with penalization (SIMP) will be used. The derivation follows the one presented in the Ph.D. Thesis of Sigmund [7].

For linear-elastic structures minimum of compliance means minimization of strain energy, which is taken as objective function. To ensure that the resulting structure has only uniaxial stress state and to simplify its engineering interpretation, some additional constraints are introduced: (i) the total amount of material (V) to be distributed in the design domain is fixed (ii) the density ($\rho_e$) can change only in certain range (bound constraints) and (iii) equilibrium constraint. After discretization with finite elements the optimization problem can be formulated as follows:

$$\begin{aligned}
& \min_{\rho,\mathbf{u}} \ \mathbf{u}^T \mathbf{K}(\rho)\mathbf{u} \\
& s.t. \\
& \mathbf{K}(\rho)\mathbf{u} = \mathbf{p} \\
& \sum_{e=1}^{N} \rho_e v_e = V \\
& 0 < \rho_e < 1
\end{aligned} \quad (1)$$

where
$\mathbf{u}$ – vector of node's displacement
$\rho$ – vector of element's densities
$\rho_e$ – density in some element
$\mathbf{K}(\rho)$ – stiffness matrix of the whole structure
$\mathbf{K}_e(\rho_e)$ – element's stiffness matrix, according to SIMP concept: $\mathbf{K}_e = (\rho_e)^p \mathbf{K}_e^{Elastic}$
   (p – penalty parameter)
$\mathbf{p}$ – vector of external loads
V – total amount of material to be distributed
$v_e$ – element's volume
e=1..N – sequence and the total number of elements

The KKT optimality condition for this problem can be formally interpreted (see Luenberger [8]) as a stationary point of the following Lagrangian function with respect to densities, displacements *and* multipliers. A *nonstrict* character of such interpretation can be illustrated by the fact that the KKT theorem states the existence of the Lagrange multipliers *only in the optimum*, at the same time, in the Lagrangian, which is defined by (2), the multipliers *do exist in each design point*.

$$L = \mathbf{u}^T \mathbf{K}\mathbf{u} + \lambda \left( \sum_{e=1}^{N} \rho_e v_e - V \right) + \boldsymbol{\mu}^T (\mathbf{K}\mathbf{u} - \mathbf{p}) + \sum_{e=1}^{N} \alpha_e (-\rho_e + 0) + \sum_{e=1}^{N} \beta_e (\rho_e - 1) \quad (2)$$

where $\lambda, \mu, \alpha, \beta$ are Lagrangian multipliers corresponding to constraints on fixed amount of material, equilibrium constraint and bound's constraints, respectively.

Using the standard procedure, which will be not repeated here (see Sigmund [7]), the optimality condition with respect to the density in each element looks as follows:

$$\begin{aligned}
& -\mathbf{u}_e^T \frac{\partial \mathbf{K}_e}{\partial \rho_e} \mathbf{u}_e + \lambda - \alpha_e + \beta_e = 0 \\
& \sum_{e=1}^{N} \rho_e v_e - V = 0 \\
& -\rho_e + 0 = 0 \ or \ \rho_e - 1 = 0
\end{aligned} \quad (3)$$

If bound constraints are not active, then the corresponding multipliers ($\alpha_e, \beta_e$) are equal to zero and the optimality condition simplifies to:

$$B_e \equiv \frac{1}{\lambda} \mathbf{u}_e^T \frac{\partial \mathbf{K}_e}{\partial \rho_e} \mathbf{u}_e = 1 \qquad (4)$$

Up to this point the derivations were mathematically well-founded. The heuristic begins by the formulating an update rule, which reads as:

$$\rho_e^{i+1} = \begin{cases} \rho_e^i B_e^i \\ \max\{(1-\zeta)\rho_e^i, 0\} & \text{if } \left(\rho_e^i B_e^i < \max\{(1-\zeta)\rho_e^i, 0\}\right) \\ \min\{(1+\zeta)\rho_e^i, 1\} & \text{if } \left(\rho_e^i B_e^i > \min\{(1+\zeta)\rho_e^i, 1\}\right) \end{cases} \qquad (5)$$

This update rule provokes some critical questions, which confirm its heuristic character:
(i) if some design point is an optimal point of the problem (1), then the design variables will not be changed by the multiplication with identity (4), but if it is not the case then it is not clear why such multiplication will lead to the optimum – all derivations were carried out under assumption of stationarity of (2);
(ii) the role of the Lagrange multipliers is also not clear and very contradictory: on the one hand, the multiplier $\lambda$, which corresponds to the fixed amount of material's constraint is calculated at the each optimization iteration during the *inner iteration loop* using Newton's method (Bendsøe & Sigmund [3], Sigmund [7]), on the other hand, the multipliers $\alpha_e$, $\beta_e$, corresponding to the bound constraints are simply ignored, if these constraints are active – the density is fixed at the bound value;
(iii) the role of the move limit parameter $\zeta$ is also not clear: it is known (from the practical experience) that it should be small, but why it is important is not explained in the classical works ([1], [3], [7])

The solution to these contradictions is the subject of the following sections.

## 2. Reformulation of OC method using projected gradient

The problem formulated in (1) belongs to the broad class of the optimization problems with state constraints in the form of partial differential equations (Gunzburger [9]). The usual way to solve such problems is to use the primal methods (Luenberger [8], Haftka [6]). The main feature of these methods is that the state constraint is enforced *at the each* optimization iteration. In contrast to them, the dual methods or "one-shot" methods (according to Gunzburger [9]) do not force the state constraints to be fulfilled at the each iteration; they try to solve the complete equation's system of KKT condition for primal variables, for state variables and for Lagrange multipliers *simultaneously*. This is done usually using the Newton's method.

Which one is suitable for explaining the working properties of optimality criteria methods can be recognized from the concrete form of KKT equation's system. The Lagrange multipliers, corresponding to the equilibrium constraint are already excluded from (3). This indicates that the interpretation must be done using primal methods.

Enforcing the state constraint means that the displacement and densities are not more independent. An *infinite* small change of the former will cause an infinite small change of the later in such a way that the equilibrium constraint remains fulfilled. This allows us to calculate their sensitivities.

$$\frac{\partial}{\partial \rho_e}(\mathbf{Ku}-\mathbf{p})=0 \Rightarrow \frac{\partial \mathbf{u}}{\partial \rho_e}=-\mathbf{K}^{-1}\left(\frac{\partial \mathbf{K}}{\partial \rho_e}\mathbf{u}\right) \tag{6}$$

The derivation of the objective function with respect to densities reads as:

$$f_{,\rho_e}=\frac{\partial}{\partial \rho_e}\left(\mathbf{u}^T\mathbf{Ku}\right)=2\mathbf{u}^T\mathbf{K}\frac{\partial \mathbf{u}}{\partial \rho_e}+\mathbf{u}^T\frac{\partial \mathbf{K}}{\partial \rho_e}\mathbf{u} \tag{7}$$

Substitution of (6) in (7) leads to the following expression for the gradient of the objective function:

$$f_{,\rho_e}=-\mathbf{u}^T\frac{\partial \mathbf{K}}{\partial \rho_e}\mathbf{u} \tag{8}$$

which is negative. If we want to *minimize* the function $f$ we have to go in descent direction, which is the negative gradient of objective function. Taking into account the local character of dependence between global stiffness matrix and density of some finite element we find the following final expression for the descent direction of the objective function:

$$-f_{,\rho_e}=\mathbf{u}_e^T\frac{\partial \mathbf{K}_e}{\partial \rho_e}\mathbf{u}_e \tag{9}$$

This expression is almost equivalent to the ones used in the optimality criteria (4). The only difference is the absence of the Lagrange multiplier $\lambda$. To understand its role we will use its alternative definition, proposed by Hestenes [10] (see also Rockafellar [11] for overview). He has introduced the Lagrange multipliers not at the optimal point of constrained optimization problem, but during *projection of the objective function's gradient onto tangential space of active constraints*.

The elementary example shown in Figure 1 illustrates this alternative definition of the multiplier $\lambda$.

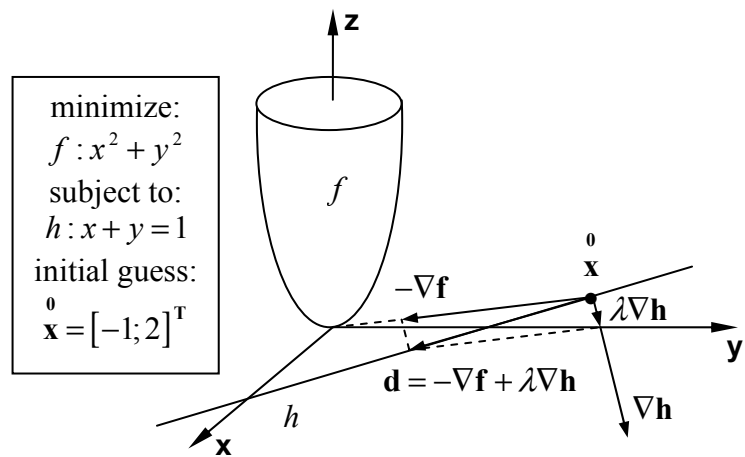

**Figure 1** Steepest descent method for a constrained problem.

To solve this problem using the steepest descent method (without line-search) a *projected gradient* has to be built. The term *"projected"* means that an *infinitely* small step in this direction will not violate the active constraint. The term *"gradient"* means that the reduction of the objective function happens in quickest possible way. These requirements are fulfilled by the following vector:

$$\mathbf{d} = -\nabla \mathbf{f} + \lambda \nabla \mathbf{h} \tag{10}$$

where $\lambda$ is a solution of the orthogonality equation:

$$(-\nabla \mathbf{f} + \lambda \nabla \mathbf{h}) \cdot \nabla \mathbf{h} = 0 \;\Rightarrow\; \lambda = \frac{\nabla \mathbf{f} \cdot \nabla \mathbf{h}}{\nabla \mathbf{h} \cdot \nabla \mathbf{h}} \tag{11}$$

It is clear that an infinitely small step along vector **d** will not violate the active constraint. In the special case, where constraint is linear, a step of *any* length (parameter $\gamma$) will not violate it:

$$\overset{i+1}{\mathbf{x}} = \overset{i}{\mathbf{x}} + \gamma \mathbf{d} \tag{12}$$

The fact that vector **d** reduces the objective function in the quickest possible way (and in this sense can be understood as gradient) will be proved in the next section.

Recalling the original problem (1) there is a simple arithmetic interpretation of parameter $\lambda$ in (10). If all components of the constraint's gradient are equal to one (for example, structured grid with elements of equal volume), then parameter $\lambda$ is equal to the mean value of the objective function's gradient. Figure 2 illustrates this.

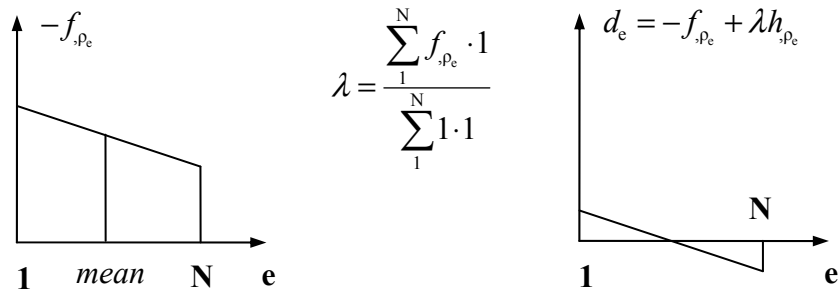

**Figure 2** Parameter $\lambda$ as a mean value of the objective function's gradient.

The OC formula (5) does actually the same job, but in an *implicit* way (in the inner iteration loop, see Sigmund [7]). It iterates until the value of Lagrange multiplier will be equal to the "weighted" mean of the strain energy (in the general case of unstructured grid, see Bendsøe & Sigmund [3]). Afterwards each density will be multiplied with the scaled derivation of objective function. If some element has a higher level of strain energy than the mean value, its density will increase ($B_e$ is higher than 1.0), if lower, then its density decreases ($B_e$ is lower than 1.0). It is really a matter of taste what to do: one can add a small positive number to the design variable or multiply it with a number higher than one.

If in some design point more than one constraint is active, the expression (10) generalizes to:

$$\mathbf{d} = -\nabla \mathbf{f} + \sum_1^S \lambda_k \nabla \mathbf{h}_k \tag{13}$$

where
k=1..S – sequence and total number of active constraints
$\nabla \mathbf{h}_k$ – gradient of active constraint
$\lambda_k$ – corresponding Lagrange (Hestenes) multiplier

The resulting vector has to be orthogonal to the *each* of the active constraint's gradient. This results in exactly S equations, which allow calculating the unique value of Hestenes multipliers. At this point it is interesting to note a full analogy to the Galerkin based Finite Element Method (Bathe [12]). The coefficients of shape functions (node's displacements) are determined in this method also using orthogonality condition between residuum of the partial differential equation, which is being solved and each of the shape functions. The inner product in this case is defined as a volume integral over the finite element. The final form of equation's system, representing orthogonality condition reads:

$$\begin{cases} \left(-\nabla \mathbf{f} + \sum_1^S \lambda_k \nabla \mathbf{h}_k\right) \cdot \nabla \mathbf{h}_1 = 0 \\ \quad\quad ... \\ \left(-\nabla \mathbf{f} + \sum_1^S \lambda_k \nabla \mathbf{h}_k\right) \cdot \nabla \mathbf{h}_s = 0 \end{cases} \Rightarrow \quad \boldsymbol{\lambda} = \left(\mathbf{HH}^T\right)^{-1}\left(\nabla \mathbf{fH}\right) \quad (14)$$

where
$\mathbf{H}$ – SxN dimensional matrix of the active constraint's gradients ($\mathbf{H}=[h_{k,\rho e}]$)
$\boldsymbol{\lambda}$ – S dimensional vector of Hestenes multipliers

For the case of simple bound constraints the orthogonality condition (14) leads to the following result: the component of the projected gradient $\mathbf{d}$, corresponding to an active constraint, must be zero. Only such a vector can be orthogonal to the gradient of the bound constraint. The following expressions illustrate this statement.

$$\nabla \mathbf{f} = \left[f_{,\rho_1}, ... f_{,\rho_{e-1}}, f_{,\rho_e}, f_{,\rho_{e+1}}, ..., f_{,\rho_n}\right]^T$$
$$h_k : \rho_k - bound = 0 \quad \Rightarrow \quad \nabla \mathbf{h}_k = [0, ..., 0, 1, 0, ...0]^T$$
$$\mathbf{d} \cdot \nabla \mathbf{h}_k = 0 \quad (15)$$
$$...0 + (f_{,\rho_{e-1}} - \lambda \cdot 0)\cdot 0 + (f_{,\rho_e} - \lambda \cdot 1)\cdot 1 + (f_{,\rho_{e+1}} - \lambda \cdot 0)\cdot 0 + 0... = 0 \quad \Rightarrow \quad \lambda = f_{,\rho_e}$$
$$\boxed{d_e = 0}$$

If we do a perturbation of design variables according to (12), than the density, where a bound constraint is active, will not change. And this is exactly what update rule (5) does.

The projected gradient method allows also to understand the importance of the move limit $\zeta$ in (5). The sensitivity analysis linearizes the objective function in an *infinitely* small area around the current design point. If we use this sensitivity for a *finite* perturbation of design variables, we have to keep it small, hoping at the same time that the linearized function is smooth enough.

## 3. Comparison with the generalized OC method

In this section it will be shown that the generalized OC method developed by Venkayya [5] is much closer to the projected gradient method than the fix-point update rule (5). Venkayya calculates Lagrangian multipliers in a way, which is almost equivalent to the orthogonality condition (14). For simplicity reasons we will show this similarity assuming three dimensional design space and two active constraints in the current design point.

The optimality condition reads as:

$$\begin{aligned} \lambda_1 \frac{\partial h_1}{\partial x} + \lambda_2 \frac{\partial h_2}{\partial x} &= \frac{\partial f}{\partial x} \\ \lambda_1 \frac{\partial h_1}{\partial y} + \lambda_2 \frac{\partial h_2}{\partial y} &= \frac{\partial f}{\partial y} \quad \Rightarrow \quad \bar{\mathbf{E}} \boldsymbol{\lambda} = \nabla \mathbf{f} \\ \lambda_1 \frac{\partial h_1}{\partial z} + \lambda_2 \frac{\partial h_2}{\partial z} &= \frac{\partial f}{\partial z} \end{aligned} \qquad (16)$$

This is the original form of the KKT condition. In Eq.(7) in [5] each equation is scaled by corresponding component of the objective function's gradient (Eq.(11) in [5]):

$$\mathbf{E}\boldsymbol{\lambda} = \mathbf{1} \qquad (17)$$

Venkayya obtains the resolving system of equations by multiplying the scaled KKT condition with the "weighted" original matrix. (see Eq.(12) in [5])

$$\mathbf{E}^T \mathbf{A} \mathbf{E} \boldsymbol{\lambda} = \mathbf{E}^T \mathbf{A} \mathbf{1} \qquad (18)$$

where **A** is a diagonal weighting matrix (Eq.(13) in [5]):

$$\mathbf{A} = \begin{bmatrix} \frac{\partial f}{\partial x} x & 0 & 0 \\ 0 & \frac{\partial f}{\partial y} y & 0 \\ 0 & 0 & \frac{\partial f}{\partial z} z \end{bmatrix} \qquad (19)$$

The expanded form of the left side of equation (18) reads

$$\mathbf{E}^T \mathbf{A} \mathbf{E} = \begin{bmatrix} \dfrac{(\nabla_x h_1)^2 x}{\nabla_x f} + \dfrac{(\nabla_y h_1)^2 y}{\nabla_y f} + \dfrac{(\nabla_z h_1)^2 z}{\nabla_z f} & \dfrac{\nabla_x h_1 \nabla_x h_2 x}{\nabla_x f} + \dfrac{\nabla_y h_1 \nabla_y h_2 y}{\nabla_y f} + \dfrac{\nabla_z h_1 \nabla_z h_2 z}{\nabla_z f} \\ \dfrac{\nabla_x h_1 \nabla_x h_2 x}{\nabla_x f} + \dfrac{\nabla_y h_1 \nabla_y h_2 y}{\nabla_y f} + \dfrac{\nabla_z h_1 \nabla_z h_2 z}{\nabla_z f} & \dfrac{(\nabla_x h_2)^2 x}{\nabla_x f} + \dfrac{(\nabla_y h_2)^2 y}{\nabla_y f} + \dfrac{(\nabla_z h_2)^2 z}{\nabla_z f} \end{bmatrix} \qquad (20)$$

The expanded form of the right side of equation (18) reads

$$\mathbf{E}^{\mathrm{T}}\mathbf{A1} = \begin{bmatrix} \dfrac{\partial h_1}{\partial x}x + \dfrac{\partial h_1}{\partial y}y + \dfrac{\partial h_1}{\partial z}z \\ \dfrac{\partial h_2}{\partial x}x + \dfrac{\partial h_2}{\partial y}y + \dfrac{\partial h_2}{\partial z}z \end{bmatrix} \quad (21)$$

A close look at expressions (18), (20) and (21) allows us to rewrite them using vector-matrix format, similar to (14), which reads:

$$\left((\mathbf{Hx})(\mathbf{H}\nabla \mathbf{f}^{-1})\right)\lambda = (\mathbf{Hx}) \quad (22)$$

The expression (22) does not allow a clear interpretation in the sense of orthogonality condition (14) yet. To simplify the interpretation we can multiply the original form (16) with a weighting matrix (19). In this case we get the expression, which states the orthogonality condition between "scaled" projected gradient and the gradients of active constraints.

$$\overline{\mathbf{E}}^{\mathrm{T}}\mathbf{A}\overline{\mathbf{E}}\lambda = \overline{\mathbf{E}}^{\mathrm{T}}\mathbf{A}\nabla \mathbf{f} \;\Rightarrow\; \left((\mathbf{Hx})\mathbf{H}\right)\lambda = \mathbf{H}(\mathbf{x}\nabla \mathbf{f})$$
$$\text{or} \quad (23)$$
$$\left(-\mathbf{x}\nabla \mathbf{f} + \mathbf{x}\lambda^{\mathrm{T}}\mathbf{H}\right)\mathbf{H} = \mathbf{0}$$

As one can see from the last equation, multiplying by the design vector **x** is actually not necessary for the definition of Lagrangian multipliers. Without it, expression (23) would be equivalent to the orthogonality condition (14).

As already mentioned, after calculating the multipliers, the design variables can be changed in two equivalent ways (depending on one's taste):
– usually used in mathematical programming (with small positive $\gamma$):

$$\overset{i+1}{\rho_e} = \overset{i}{\rho_e} + \gamma\left(-f_{,\rho_e} + \sum_1^S \lambda_k h_{k,\rho_e}\right) \quad (24)$$

– usually used in optimality criteria methods (with small positive $\gamma$):

$$\overset{i+1}{\rho_e} = \overset{i}{\rho_e}\left(\sum_1^S \lambda_k \frac{h_{k,\rho_e}}{f_{,\rho_e}}\right)^{\gamma} \quad (25)$$

The second basic element of the OC algorithm according to [5] is the scaling of new design vector in order to fulfill the violated constraints. Its necessity becomes clear if we remember that we do a *finite* step in tangential direction of the constraints. The only exception where scaling is not necessary would be linear constraints.

## 4.    Least squares estimates or *real* multipliers?

The expression (14) is actually well known in the optimization theory (see Haftka [6]), but it is understood in a different way – as *least-square estimate* of Lagrange multipliers. For example,

Fletcher has proposed in [13] the following augmented Lagrangian to solve constrained optimization problems:

$$L^F = -f + \sum_1^S \lambda_k h_k + \frac{1}{2\mu}\sum_1^S h_k h_k \qquad (26)$$

where the multipliers $\lambda_k$ minimize the square norm of the vector, which we call projected gradient. Indeed, the expression (14) can be obtained as a solution of the following auxiliary problem:

$$\min_\lambda \frac{1}{2}(-\nabla \mathbf{f} + \lambda \mathbf{H})\cdot(-\nabla \mathbf{f} + \lambda \mathbf{H}) \;\Rightarrow\; \frac{\partial F}{\partial \lambda} = (-\nabla \mathbf{f} + \lambda \mathbf{H})\mathbf{H} = 0$$
$$\lambda = (\mathbf{H}\mathbf{H}^T)^{-1}(\mathbf{H}\nabla \mathbf{f}) \qquad (27)$$

Formally, it is not wrong to understand the Hestenes multipliers as least-square estimates of Lagrange multipliers. The classical KKT theorem introduces the multipliers *only in the optimum*, so any expression for their calculation in the non-optimal point can be only an "estimation".

In contrast to that, the main idea of this paper is that the Hestenes multipliers from (14) are not estimates, but the *real* multipliers. And only in special case if resulting projected gradient is zero vector, they obtain the meaning of classical Lagrange multipliers. Such interpretation obviously requires reformulation of Karush-Kuhn-Tucker condition, which will be done in the following section.

## 5. Alternative formulation of KKT condition

The projected gradient resulting from (14) can be directly used for an alternative definition of the optimal point of optimization problem with constraints. This condition takes the following simple form:

$$\mathbf{d} = \mathbf{0} \qquad (28)$$

This reformulation has the following advantages:
(i) full analogy with unconstrained optimization (gradient is a zero vector);
(ii) multipliers exist not only in optimum, but in each design point. They also have a clear geometric interpretation (see Figure 1);
(iii) the projected gradient can be directly used in some "steepest descent like" optimization algorithms. It was shown that these algorithms are equivalent to the optimality criteria methods.

It is important to prove that there is no other vector, which reduces the objective function without violation of constraints. It can be proven by contradiction, assuming that there is some vector **s**, which is linear independent on **d**. Such a vector can be always constructed by the following linear combination:

$$\mathbf{s} = \alpha\,\mathbf{d} + \beta\,\mathbf{r} \qquad (29)$$

where the vectors **d** and **r** are orthogonal, i.e. **d**•**r**=0, or

$$\left(-\nabla \mathbf{f}+\sum_{1}^{S}\lambda_{k}\nabla \mathbf{h}_{k}\right)\cdot \mathbf{r}=0 \qquad (30)$$

At the same time, vector **s** must not violate any of active constraints. This means it must be orthogonal to the each constraint's gradients, i.e. $\mathbf{s}\cdot\nabla\mathbf{h}_k=0$ or

$$(\alpha\mathbf{d}+\beta\mathbf{r})\cdot\nabla\mathbf{h}_k=0 \quad \Rightarrow \quad \alpha\mathbf{d}\cdot\nabla\mathbf{h}_k+\beta\mathbf{r}\cdot\nabla\mathbf{h}_k=0 \qquad (31)$$

but according to (14) we have $\mathbf{d}\cdot\nabla\mathbf{h}_k=0$ from which follows that $\mathbf{r}\cdot\nabla\mathbf{h}_k=0$. Substituting the last result in (30) one obtains that the vector **r** is orthogonal to the objective function's gradient.

$$\nabla\mathbf{f}\cdot\mathbf{r}=0 \qquad (32)$$

This orthogonality proofs that there is no other vector than **d**, which reduces the objective function without violation of constraints. Therefore, if this vector is a zero vector, then the current design point is an optimal point of constrained optimization problem.

## 6.   Design of tension-only structures

To verify the developed method it will be applied to one engineering problem: topology design of tension-only structures. Its practical application is conceptual design of cable supported bridges. The problem was already considered by Meiss in [14], where the optimality criteria method was used. Here we present a more rigorous way (using numerical integration) to suppress the compressive stress state.

Recalling again the original problem (1) and taking into account the SIMP-type dependence between Young's modulus and density of material in some element, it follows from (9) that the derivative of the strain energy, stored in the element, is proportional to itself:

$$-f_{,\rho_e}=\mathbf{u}_e^T\frac{\partial \mathbf{K}_e}{\partial \rho_e}\mathbf{u}_e=\frac{p}{\rho_e}\mathbf{u}_e^T\mathbf{K}_e\mathbf{u}_e=\frac{p}{\rho_e}\int_V \mathbf{u}_e^T\mathbf{B}^T\mathbf{C}\mathbf{B}\mathbf{u}_e\,\det\mathbf{J}\,dV \qquad (33)$$

The integration of the stiffness matrix is usually done numerically, i.e. the integral in (33) is replaced by Gauss quadrature. For example, for a 2D-bilinear element the four-point integration scheme is exact (see Bathe [12]).

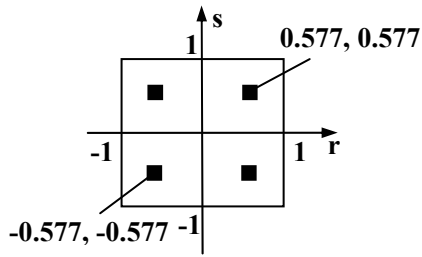

$$\int_\Omega \mathbf{u}_e^T\mathbf{B}^T\mathbf{C}\mathbf{B}\mathbf{u}_e\,\det\mathbf{J}\,d\Omega=\sum_2\sum_2 1\cdot\mathbf{u}_e^T\mathbf{B}^T\mathbf{C}\mathbf{B}\mathbf{u}_e\,\det\mathbf{J} \qquad (34)$$

At the same time, the stored elastic energy can be calculated not only through strains, but also through stresses:

$$\int_V \mathbf{u}_e^T \mathbf{B}^T \mathbf{C B u}_e \det \mathbf{J} dV = \int_V \mathbf{\varepsilon} : \mathbf{C} : \mathbf{\varepsilon} \det \mathbf{J} dV = \int_V \mathbf{\sigma} : \mathbf{C}^{-1} : \mathbf{\sigma} \det \mathbf{J} dV \qquad (35)$$

Numerical integration of the last integral in (35) requires evaluation of stresses only at the Gauss points. After calculating the stress tensor, one can perform principal axes transformation. The value of elastic energy will not change of course, but after this transformation we have a clear distinction between tensile and compressive part of the stored energy. So we can easily suppress unwanted stress state. The final expression for the descent direction (for 2D case) reads:

$$-f_{,\rho_e} = \frac{p}{\rho_e}\int_{V_e} \mathbf{\sigma} : \mathbf{C}^{-1} : \mathbf{\sigma} \det \mathbf{J} dV_e = \frac{p}{\rho_e}\sum_2\sum_2 \frac{1}{\rho_e^p E_0}\cdot\left(\overline{\sigma}_I^2 - 2\mu\overline{\sigma}_I\overline{\sigma}_{II} + \overline{\sigma}_{II}^2\right)\det \mathbf{J} \qquad (36)$$

where $\overline{\sigma}_I$ and $\overline{\sigma}_{II}$ are reduced principal stresses, according to the following scheme.

$$\overline{\sigma}_i = \begin{cases} \sigma_i & if \quad \sigma_i > 0 \\ k\sigma_i & else \end{cases}, \quad \begin{array}{l} i=I,II \\ k \in [0,1) \end{array} \qquad (37)$$

where $k$ is a reduction factor

Using this modified objective function's gradient (36) the projected gradient according to (13) and (14) is built. If boundary and loading conditions principally allow tension-only structure, then algorithm converges quickly. The Figure 3 shows one example of optimal design using this method.

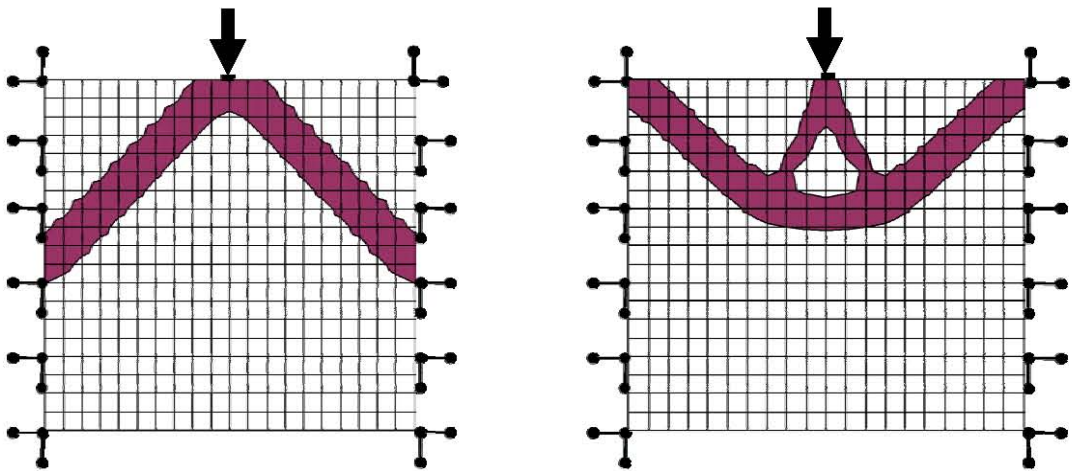

**Figure 3** Example of the standard topology and the topology with suppressed compressive stresses.

## 7. Closure

The motivation at the beginning of this work was to understand, why optimality criteria methods work. But during research it became clear that the main reason of why these methods are considered to be heuristic is the implicit character of Karush-Kuhn-Tucker condition. This theorem states the existence of Lagrangian multipliers only in the optimum. The alternative formulation of optimality condition (28) using Hestenes multipliers removes this implicit character and simplifies the understanding of constrained optimization theory. This can be probably considered as the main result of the presented work.

**Acknowledgements**

The financial support to the author from the German Academic Exchange Service (DAAD) is gratefully acknowledged. The author is also grateful to the anonymous reviewers for the critical comments and suggestions.

**References**

[1] Bendsøe M.P., Kikuchi N. 1988. Generating optimal topologies in structural design using a homogenization method. *Computer Methods in Applied Mechanics and Engineering*, 71: 197-224.

[2] Soto C.A. 2002. Applications of Structural Topology Optimization in the Automotive Industry: Past, Present and Future. In: Mang H.A., Rammerstorfer F.G., Eberhardsteiner J. *Fifth World Congress on Computational Mechanics*, Vienna, Austria.

[3] Bendsøe M.P., Sigmund O. 2003. *Topology Optimization. Theory, Methods and Applications*. Springer, Berlin.

[4] Maute K. 1998. *Optimization of Topology and Form of slender structures*. Ph.D. Thesis. University of Stuttgart. (in German)

[5] Venkayya V.B. 1993. Generalized optimality criteria method. In: Kamat M. (ed.). *Structural optimization: status and promise*. American Institute of Aeronautics and Astronautics.

[6] Haftka R.T. 1992. *Elements of Structural Optimization*. Kluwer Academic Publishers, Dordrecht.

[7] Sigmund O. 1994. *Design of materials and structures using topology optimization*. Ph.D. Thesis, TU Denmark.

[8] Luenberger D. 1989. *Linear and nonlinear programming*. Addison-Wesley, Second edition.

[9] Gunzburger M.D. 2003. *Perspectives in Flow Control and Optimization*. SIAM, Philadelphia.

[10] Hestenes M.R. 1975. *Optimization theory. The finite dimensional case*. John Wiley & Sons, New York.

[11] Rockafellar R.T. 1993. Lagrange multipliers and optimality. *SIAM Review*, 35(2): 183-238.

[12] Bathe K.-J. 1996. *Finite element procedures*. Englewood Cliffs, NJ: Prentice Hall.

[13] Fletcher R. 1973. An exact penalty function for nonlinear programming with inequalities. *Mathematical Programming*, 5: 129-150.

[14] Meiss K. 2002. New developments in engineering application of structural optimization. In: Schießl P. (ed.). *4th International Ph.D. Symposium*. Springer-VDI, Düsseldorf.